\documentclass[10pt]{article}
\usepackage[utf8]{inputenc}
\usepackage{latexsym,amssymb,amsmath,url,authblk}

\def\N{\mathbb{N}}
\def\Q{\mathbb{Q}}
\def\Z{\mathbb{Z}}
\def\R{\mathbb{R}}
\def\C{\mathbb{C}}

\def\proof{\par\noindent{\em Proof. }}
\def\eproof{\hfill{$\Box$}\bigskip}

\def\ds{\dots}
\def\sus{\subset}
\def\al{\alpha}
\def\be{\beta}
\def\ga{\gamma}

\def\cc{\colon}

\newtheorem{thm}{Theorem}[section]
\newtheorem{prop}[thm]{Proposition}

\newtheorem{defi}[thm]{Definition}

\title{The transcendence of $\mathrm{e}$ via formal power series}
\author{Martin Klazar\footnote{{\tt klazar@kam.mff.cuni.cz}}}
\date{}

\begin{document}

\maketitle

\begin{abstract}
We review Hilbert's classical analytical
proof of the transcendence of the number $\mathrm{e}$. Then, we show how this result can be obtained algebraically by means of formal power 
series (FPS). We give two proofs of the transcendence of $\mathrm{e}$ based on FPS. 
The first of them is a~specialization of the 1990 proof by 
Beukers, B\'ezivin and Robba of the Lindemann--Weierstrass theorem. The second proof is due to this author and is 
an adaptation of Hilbert's argument to FPS.    
\end{abstract}

\section{Hilbert's proof and how to make it countable}\label{sec_intro}

The transcendence of the number $\mathrm{e}$ was first proven by Hermite \cite{herm}.
Hilbert's proof \cite{hilb} goes as 
follows. We suppose for the contradiction that there exists 
$n+1$ integers $a_0$, $a_1$, $\ds$, $a_n$, not all of 
them $0$, such that  
$$
P:=
a_0+a_1\mathrm{e}+a_2\mathrm{e}^2+\ds
+a_n\mathrm{e}^n=0\,.
$$
We may assume that $a_0\ne0$. For any $r\in\N$ ($=\{1,2,\ds\}$) we define integral polynomials 
$$p_r(x)=
x^r((x-1)(x-2)\ds(x-n))^{r+1}\,, 
$$
with $n$ as in $P$. We have the expansions 
$$
p_r(x)=(-1)^{n(r+1)}(n!)^{r+1}\cdot x^r+ax^{r+1}+\ds
$$
and
$$
p_r(x+j)=bx^{r+1}+cx^{r+2}+\ds\,,
$$
where $j=1,2,\ds,n$ and $a,b,c,\ds$ are in $\Z$ (the integers). We consider improper integrals
$${\textstyle
I_r=\int_0^{+\infty}p_r(x)\cdot\mathrm{e}^{-x}\,\mathrm{d}x\,.
}
$$
Using their additivity,  we express
the product $PI_r$ as
\begin{eqnarray*}
0&=&PI_r=A_r+B_r\\
&:=&{\textstyle
\sum_{j=0}^n a_j\mathrm{e}^j\int_0^j 
p_r(x)\cdot\mathrm{e}^{-x}\,\mathrm{d}x
+\sum_{j=0}^n
a_j\mathrm{e}^j\int_j^{+\infty}
p_r(x)\cdot\mathrm{e}^{-x}\,\mathrm{d}x
}\,.
\end{eqnarray*}
As for the quantity $A_r=\sum_{j=0}^n a_j\mathrm{e}^j\int_0^j 
p_r(x)\cdot\mathrm{e}^{-x}\,\mathrm{d}x$, it is easy to 
estimate it as 
$$
|A_r|\le c^r,\ r\in\N\,,
$$
where $c\ge1$ is a~constant depending only on the coefficients $a_0$, $a_1$, $\ds$, $a_n$. 

We claim that the second quantity
$$
{\textstyle
B_r=\sum_{j=0}^n
a_j\mathrm{e}^j\int_j^{+\infty}
p_r(x)\cdot\mathrm{e}^{-x}\,\mathrm{d}x
}
$$
is an integral multiple of $r!$ and that $B_r\ne0$ 
for infinitely many $r\in\N$. Since $c^r/r!\to0$ as $r\to\infty$ for any 
constant $c\ge1$, the above equalities
$$
A_r+B_r=0\ \ (r\in\N)
$$
are contradictory and we are done. We prove the claim by means of the integral identity  
$${\textstyle
\int_0^{+\infty}x^n\mathrm{e}^{-x}\,\mathrm{d}x=n!\ \ 
(n\in\N_0=\{0,1,\ds\})\,.
}
$$ 
It is due to Euler and is easy to 
establish by induction on $n$ via integration by parts. The linearity of improper integrals yields the more general identity 
$$
{\textstyle
\int_0^{+\infty}\big(\sum_{j=0}^m b_jx^j\big)\mathrm{e}^{-x}\,\mathrm{d}x=\sum_{j=0}^m b_j\cdot j!\ \ 
(m\in\N_0,\,b_j\in\Z)\,.
}
$$
We compute that
\begin{eqnarray*}
B_r&=&{\textstyle
\sum_{j=0}^n
a_j\mathrm{e}^j\int_j^{+\infty}
p_r(x)\cdot\mathrm{e}^{-x}\,\mathrm{d}x
}\\
&=&{\textstyle
\sum_{j=0}^n
a_j\int_j^{+\infty}
p_r((x-j)+j)\cdot\mathrm{e}^{-(x-j)}\,\mathrm{d}x
}\\
&=&{\textstyle
\sum_{j=0}^n
a_j\int_0^{+\infty}
p_r(y+j)\cdot\mathrm{e}^{-y}\,\mathrm{d}y}\\
&=&\pm a_0(n!)^{r+1}\cdot r!+m_r\cdot(r+1)!\ \ (m_r\in\Z)\,.
\end{eqnarray*}
The first equality follows from the definition of $B_r$.
In the second equality, we use the exponential identity 
$\mathrm{e}^{a+b}=\mathrm{e}^a\mathrm{e}^b$, the linearity of improper integrals, and an algebraic rearrangement. In the third 
equality, we change the variable $x$ to $y=x-j$.  
In the fourth equality, we use the above generalization of Euler's identity and
the expansions of $p_r(y+j)$ for $j=0,1,\ds,n$. Clearly,  
$B_r\in\Z$ and is divisible by $r!$. If $B_r=0$ then 
$r+1$ divides $a_0(n!)^{r+1}$. This divisibility is excluded when $r+1$ is 
coprime to the number $a_0\cdot n!$. Hence $B_r\ne0$ for infinitely many 
$r\in\N$. 
\eproof

\noindent
Baker \cite[Theorem~1.2]{bake} 
gives a~variant of this proof that uses only proper integrals 
$\int_a^b$. A~similar proof appears in Conrad \cite{conr}. In \cite{hilb}, Hilbert proved in a~similar way also the 
transcendence of $\pi$.

Most research mathematicians do not care about using or not using uncountable sets in their work, but for the minority of 
those who do, Hilbert's proof poses the natural challenge of eliminating all substantial uses of uncountable sets from it. These are uses of the functions 
$$
F(x)=p(x)\cdot\mathrm{e}^{-x}\cc[0,\,+\infty)\to\R\ \ (p(x)\in\Z[x])\,.
$$  
Four explanations are in order. 1. It should be intuitively clear what is 
a~substantial use of an uncountable 
set. For example, in the result that $\sum_{n\ge0}x^{-n}=\frac{1}{1-x}$ for every $x\in\R$ with $|x|<1$, the mention of the set 
of real numbers $\R$ is 
{\em not} a substantial use of an uncountable set because we can easily formulate it away. 2. In mathematical logic, it has 
been known probably since the publication of \cite{weyl} that real analysis is a~countable undertaking 
because it can be carried out in systems of second order arithmetic. 
In our elimination of uncountable sets from Hilbert's proof, we are not interested in approaches of mathematical logic. We are 
interested  in technically simpler approaches where one does 
not need to know, for example, what 
a~formula is. 3. After the elimination of substantial uses of uncountable sets 
from the proof, the remaining substantially used sets should not only be at most countable; they 
are required to be hereditarily at most countable. 4. In order to meet the last condition, individual real numbers 
are represented by Dedekind cuts.

What would Hilbert suggest to us if we could ask him? Perhaps: ``If you do not like uncountable sets in my proof, which I~respect because I myself devoted considerable effort to eliminating infinity from proofs 
(\cite{zach}), you can get rid of them easily. Those functions $F(x)$ are quite special; they are analytic. If you 
need a~functional value $F(b)$ at a~real number $b\ge0$, you do not have to search through the uncountable set $F$ to pick up 
the unique pair $\langle b,F(b)\rangle\in F$. You can take the sequence of coefficients $(a_0,a_1,\ds)$ in the 
expansion $F(x)=\sum_{n\ge0}a_nx^n$, which is a~(hereditarily at most) countable set, and compute the limit $\lim_{n\to\infty}\sum_{j=0}^n 
a_jb^j=F(b)$. Just replace $F$ with $(a_0,a_1,\ds)$.'' 

We implement Hilbert's suggestion in Section~\ref{sec_klazar}, where we adapt Hilbert's proof to the algebra 
$\R\{x\}$ of convergent FPS with real coefficients. The technical details are not as 
straightforward as one might think. When preparing our article, we found that a~countable FPS proof of the transcendence of 
$\mathrm{e}$, in fact, of a~much more general result, has been known for more than three decades. In 1990, Beukers, B\'ezivin 
and Robba \cite{beuk_al} devised an FPS proof of the Lindemann--Weierstrass theorem (for the precise statement of it, see the next 
section). Since their proof is little known, we review its specialization to the transcendence of~$\mathrm{e}$ in 
Section~\ref{sec_beukers}. The motivation for the proof in \cite{beuk_al} was not the elimination of uncountable sets; at 
least nothing to this effect is said in \cite{beuk_al}. The 
main idea of the proof is completely different from Hilbert's: if the 
number $\mathrm{e}$ were algebraic, a~differential equation would have 
a~rational FPS solution $V(x)$,  which can be shown to be impossible 
by considering the poles. 

Another approach for obtaining countable proof of the transcendence of the number $\mathrm{e}$, not based on FPS, is implemented 
in \cite{klaz}. The idea is simple: replace 
$F\cc[0,+\infty)\to\R$ with the (hereditarily at most) countable restriction 
$F\,|\,[0,+\infty)\cap\Q$. In \cite{klaz}, we more generally develop real analysis for the {\em countable} real functions 
$f$ of the form $f\cc M\to\R$ where $M\sus\Q$. One may see some problems with this approach: for example, it is easy to 
define a~continuous function $f\cc[0,1]\cap\Q\to\R$ that is unbounded. So the minimax principle fails for 
countable functions. Yes, but we only work with (locally) uniformly continuous (UC) countable functions, and the 
minimax principle is easily adapted to them. The real problem is different: to adapt the standard result to (locally) UC  countable functions, that if $b\in X\sus\R$, $b$ 
is a~two-sided limit point of $X$, and $f\cc X\to\R$ is such that
$f'(b)\ne0$, then the value $f(b)$ is not a~local extreme of the function 
$f$.  See \cite{klaz} for details.

\section{Specializing the proof of Beukers, B\'ezivin and Robba in \cite{beuk_al}}\label{sec_beukers}

The Lindemann--Weierstrass (LW) theorem says that
$$
b_1\mathrm{e}^{\al_1}+b_2
\mathrm{e}^{\al_2}+\ds+b_t\mathrm{e}^{\al_t}\ne0
$$
whenever $b_j,\al_j\in\C$ are $2t\ge2$ algebraic numbers such that 
$b_j\ne0$ and the $\al_j$ are mutually distinct. See 
Baker \cite[Theorem~1.4]{bake} for an analytic proof. Using 
formal power series, Beukers, B\'ezivin, and Robba \cite{beuk_al} obtained another proof of the LW theorem. In this section, we 
specialize their proof to $b_j,\al_j\in\Z$, which 
gives the transcendence of $\mathrm{e}$. In order to clearly 
present the logical structure of the proof, we reorganize 
the exposition in \cite{beuk_al}. 

We first recall that a~FPS $f(x)\in \C[[x]]$ is rational if there exist 
two polynomials $p(x),q(x)\in\C[x]$ such that $q(0)=1$ and $q(x)f(x)=p(x)$. We 
may assume that $p(x)$ and $q(x)$ are coprime. If $f(x),g(x)\in\C[[x]]$ are 
rational, then so are $f(x)+g(x)$, $f(x)g(x)$, and $f'(x)$. Any rational 
FPS $f(x)\in\C[[x]]$ has a~unique decomposition into partial fractions
$${\textstyle
f(x)=r(x)+\sum_{j=1}^t\sum_{i=1}^{m_j}
\frac{\be_{j,i}}{(1-\al_jx)^i}\,,
}
$$
where $r(x)\in\C[x]$, $t\in\N_0$, $m_j\in\N$, $\be_{j,i},\al_j\in\C$ with $\be_{j,m_j}\ne0$, and 
the $\al_j$ are nonzero and mutually distinct. For $t=0$, the double sum is defined as $0$. The numbers $t$, $m_j$, and $\al_j$ come from the factorization
$$
{\textstyle
q(x)=\prod_{j=1}^t(1-\al_jx)^{m_j}\,.
}
$$
The partial fractions are the formal power series
$${\textstyle
\frac{\be_{j,i}}{(1-\al_jx)^i}=
\sum_{n\ge0}\binom{-i}{n}\be_{j,i}(\al_j)^n\cdot x^n\,.
}
$$
We say that $1/\al_j$ is a~pole of $f(x)$ with order $m_j$. 

\begin{prop}\label{prop_oPolech}
If $p(x)\in\C[x]$, $m\in\N_0$, $c\in\C$ with $c\ne0$, and if $A(x)$ in $\C[[x]]$ is rational, then the
rational {\em FPS}
$$
p(x)A(x)+cx^mA'(x)
$$
is either a~polynomial and has no poles, or all its poles have orders 
at least $2$.
\end{prop}
\proof
The product 
$p(x)A(x)$ only possibly decreases the orders of the poles of $A(x)$ 
or cancels them and does not create any new poles.  The product and derivative  
$cx^mA'(x)$ preserve all poles of $A(x)$, increase their orders by 
$1$, and do not create any new poles. The sum $p(x)A(x)+cx^mA'(x)$ merges the 
two sets of poles. It follows that the resulting set of poles with their 
orders is that of $cx^mA'(x)$.
\eproof

We extracted the next theorem on rational FPS from the 
arguments in \cite{beuk_al}. 

\begin{thm}\label{thm_rational}
Let $t\in\N$, $b_1$, $b_2$, $\ds$, $b_t$ be nonzero integers and $\al_1$, $\al_2$, $\ds$, $\al_t$ be distinct natural numbers. Let $n\in\N_0$. We define
integers
$$
u_n=b_1\al_1^n+b_2\al_2^n+\ds+b_t\al_t^n\,\text{ and }\,v_n=n!\cdot{\textstyle
\sum_{r=0}^n u_r/r!}\,.
$$
Suppose that for some constant $A\ge1$ we have bound
$$
v_n=O(A^n)\,\text{ for }\,n\in\N_0\,. 
$$
Then the {\em FPS}   $V(x)=\sum_{n\ge0}v_nx^n$ is rational,  because 
$$
{\textstyle
\big(1-a_1x-a_2x^2-\ds-a_tx^t\big)^k\sum_{n\ge0}v_nx^n\in\Z[x]
}
$$
for some $k\in\N$, where the integers $a_j$ are the coefficients in
$$
{\textstyle
\prod_{j=1}^t(1-\al_jx)=1-a_1x-a_2x^2-\ds-a_tx^t\,.
}
$$
\end{thm}
\proof
Let $k,n\in\N_0$. We define integers $v_n(k)$ by the relations
$$
{\textstyle
\sum_{n\ge0}v_n(k)x^n=
\big(1-a_1x-a_2x^2-\ds-a_tx^t\big)^k\sum_{n\ge0}v_nx^n\,.
}
$$
Then for every $n\ge t$ and $k\in\N_0$ we have the recurrence
$$
v_n(k+1)=v_n(k)-a_1v_{n-1}(k)-\ds-a_tv_{n-t}(k)\,.
$$
Let $k,n\in\N_0$ with $n\ge tk$. We establish two properties of the numbers $v_n(k)$.
\begin{enumerate}
\item We have the bound $|v_n(k)|\le cA^nC^k$ where $C=1+|a_1|+\ds+|a_t|$ and $c\ge0$ is a~constant.
\item The number $k!$ divides $v_n(k)$.
\end{enumerate}
The bound in part~1 follows by induction on $k$
from the assumed bound on $v_n=v_n(0)$. For $k=0$ the bound 
holds because it is the assumed bound on $v_n$. Let $k,n\in\N_0$ with $n\ge 
t(k+1)$. Then $n-t\ge tk$ and
by the above recurrence and induction, 
\begin{eqnarray*}
|v_n(k+1)|&=&|v_n(k)-a_1v_{n-1}(k)-\ds-a_tv_{n-t}(k)|\\
&\le&c\big(A^nC^k+A^{n-1}C^k|a_1|+\ds+
A^{n-t}C^k|a_t|\big)\le cA^nC^{k+1}\,.
\end{eqnarray*}

In order to prove the  divisibility in part~2, we note that $v_0=u_0$ and that for $n\in\N$,
$$
{\textstyle
\frac{v_n}{n!}-\frac{v_{n-1}}{(n-1)!}=
\frac{u_n}{n!},\,\text{ so that $v_n-nv_{n-1}=u_n$}\,.
}
$$
For $k=0$ the divisibility trivially holds. Let $k\in\N$. We 
set $(0)_0=(n)_0=1$ and $(n)_m=n(n-1)\ds(n-m+1)$ for $m\in\N$. For $n\in\N_0$ we 
combine $\min(n,k)$ displayed differences, multiply the result by $n!$, and get that
$${\textstyle
v_n=(n)_0u_n+(n)_1u_{n-1}+\ds+(n)_{k-1}u_{n-k+1}+(n)_kv_{n-k}}\,,
$$
where the $u_{\ds}$ and $v_{\ds}$ with negative indices are defined as 
$0$. We denote the last summand $(n)_kv_{n-k}$ ($\in\Z$) by $w_n$. The number $k!$ always divides $w_n$ because 
$w_n=k!\cdot\binom{n}{k}v_{n-k}$ for $n\ge k$, and $w_n=0$ for $0\le n<k$. The difference of the FPS 
$V(x)=\sum_{n\ge0}v_nx^n$ and $W(x)=\sum_{n\ge0}w_nx^n$ is
$${\textstyle
V(x)-W(x)=\sum_{n\ge0}\big((n)_0u_n+(n)_1u_{n-1}+\ds+(n)_{k-1}u_{n-k+1}\big)x^n\,.
}
$$
We claim that
$$
{\textstyle
V(x)-W(x)=\frac{p(x)}{(1-a_1x-\ds-a_tx^t)^k}=:\frac{p(x)}{q(x)}
}
$$
for some polynomial $p(x)\in\Z[x]$ with degree less than $tk$. Then for 
every $n\ge tk$ we have
$$
v_n(k)=[x^n]\,q(x)V(x)=[x^n]\,q(x)W(x)\,,
$$
where $[x^n]f(x)$ denotes the coefficient of $x^n$ in $f(x)$, 
and $k!$ divides $v_n(k)$ because every coefficient in $W(x)$ is divisible by $k!$. We prove the claim: for every $r$
with $0\le r\le k-1$ we indeed have,
since $\binom{n}{r}=(-1)^{n-r}\binom{-r-1}{n-r}$ for $n\ge r$, that
\begin{eqnarray*}
{\textstyle\sum_{n\ge0}(n)_ru_{n-r}x^n}&=&{\textstyle
\sum_{j=1}^t 
r!\cdot b_jx^r\sum_{n\ge r}
\binom{n}{r}\al_j^{n-r}x^{n-r}}\\
&=&{\textstyle
\sum_{j=1}^t r!\cdot b_jx^r
(1-\al_j x)^{-r-1}=
\frac{p_r(x)}{(1-a_1x-\ds-a_tx^t)^{r+1}}
}
\end{eqnarray*}
for some polynomial $p_r(x)\in\Z[x]$ with degree less than $t(r+1)$.

Now the two properties of numbers $v_n(k)$ imply that if $n\ge tk$ and 
$v_n(k)\ne0$, then
$$
k!\le|v_n(k)|\le cA^nC^k\,.
$$
In other words, if $n\ge tk$ and if
$k!>cA^nC^k$, then $v_n(k)=0$. Using this implication and the modification
$$
v_n(k)=v_n(k+1)+a_1v_{n-1}(k)+\ds+a_tv_{n-t}(k),\ n\ge t\,,
$$
of the above recurrence, we show that there exists a~$k\in\N$ such that 
$v_n(k)=0$ for every large $n$.
This will prove that the formal power series $V(x)=\sum_{n\ge0}v_nx^n$ is 
rational. So let $k_0\in\N$ be such that $k!>cA^{2tk}C^k$ for every $k\ge k_0$, and let
\begin{eqnarray*}
N&=&\{\langle k,\,n\rangle\in\N^2\cc\;
k\ge k_0\;\&\;tk\le n\le 2tk\}\,\text{ and}\\
M&=&\{\langle k,\,n\rangle\in\N^2\cc\;k\ge k_0\;\&\;n\ge 2tk\}\,.    
\end{eqnarray*}
Thus $\langle k,n\rangle\in N$ $\Rightarrow$ $v_n(k)=0$. Let $\langle k,n\rangle\in\N^2$. We consider the norm
$\|\langle k,n\rangle\|=n-2tk$ ($\in\Z$). We name the pairs in the modified recurrence:
$$
\ga=\langle k,\,n\rangle,\,
\be_0=\langle k+1,\,n\rangle,\,
\be_1=\langle k,\,n-1\rangle,\,\ds,\,
\be_t=\langle k,\,n-t\rangle\,.
$$
The following three properties of the norm $\|\ds\|$ and 
the pairs $\ga$, $\be_0$, $\ds$, $\be_t$ are easy to prove.
\begin{enumerate}
\item The norm $\|\ds\|$ attains a~minimum value on the set $M$ (namely, zero). 
\item If $\ga\in M$ then $\be_i\in M\cup N$ for every $i=0,1,\ds,t$.
\item For every $i=0,1,\ds,t$, the norm $\|\be_i\|<\|\ga\|$. 
\end{enumerate}
From these three properties, from the implication for the set $N$ and from the modified recurrence, it follows by 
induction on the norm that $v_n(k)=0$ for every $\langle k,n\rangle\in M$. 
This concludes the proof because we see that every number $k\ge k_0$ has the 
required property.
\eproof

We prove, with the help of Proposition~\ref{prop_oPolech} and 
Theorem~\ref{thm_rational}, that the number~$\mathrm{e}$ is transcendental. 
We assume, for the contradiction, that $b_j,\al_j\in\Z$, $j=1,2,\ds,t$,
are $2t\ge2$ integers such that
$b_j\ne0$, the $\al_j$ are mutually distinct, and that 
$${\textstyle
\sum_{j=1}^t b_j\mathrm{e}^{\al_j}=0\,.
}
$$
We may assume that $\al_j\in\N$. We set $A:=\max(\{\al_1,\ds,\al_t\})$ and define $u_n$ and $v_n$ as in Theorem~\ref{thm_rational}.  Then
$$
{\textstyle
\sum_{n\ge0}\frac{u_n}{n!}x^n=
\sum_{j=1}^t b_j\sum_{n\ge0}\frac{\al_j^n}{n!}x^n
=\sum_{j=1}^t b_j\mathrm{e}^{\al_jx}\ \ (\in\C[[x]])
}
$$
and
$$
{\textstyle
\sum_{n\ge0}\frac{u_n}{n!}=
\sum_{j=1}^t
b_j\mathrm{e}^{\al_j}=0\,.
}
$$
From the definition of $u_n$ and $A$ it easily follows that $u_n=O(A^n)$ for $n\in\N_0$.
Since by the definition of $v_n$ we have
$$
|v_n|={\textstyle
n!\cdot\big|\sum_{r=0}^n u_r/r!\big|=n!\cdot\big|\sum_{r=n+1}^{\infty}u_r/r!\big|}
$$
---\,we invoke the contradictory assumption\,---\,we get the 
bound 
\begin{eqnarray*}
|v_n|&\le&{\textstyle\frac{|u_{n+1}|}{n+1}+
\frac{|u_{n+2}|}{(n+1)(n+2)}+\ds\le
cA^n\big(
\frac{A}{n+1}+\frac{A^2}{(n+1)(n+2)}+\ds\big)}\\
&\le&{\textstyle cA^n\sum_{j\ge0}\frac{A^j}{j!}=
c\,\mathrm{e}^A\cdot A^n\,,
}
\end{eqnarray*}
where $c>0$ is an absolute constant.
Thus $v_n=O(A^n)$ for $n\in\N_0$.
Let $V(x)=\sum_{n\ge0}v_nx^n$ ($\in\Z[[x]]$). We know from the proof  of Theorem~\ref{thm_rational} that $v_n-nv_{n-1}=u_n$ for $n\in\N$ and $v_0=u_0$. In terms of FPS we get
$${\textstyle
(1-x)V(x)-x^2V'(x)=\sum_{n\ge0}(v_n-nv_{n-1})x^n=\sum_{n\ge0}u_nx^n=
\sum_{j=1}^t\frac{b_j}{1-\al_jx}
}
$$
(with $v_{-1}:=0$), which is a~rational FPS with $t$ poles of 
order $1$ each. This
contradicts Proposition~\ref{prop_oPolech} 
because $V(x)$ 
is rational by Theorem~\ref{thm_rational}.
\eproof

\noindent
It is clear that the previous proof does not substantially use 
uncountable sets. In the original proof of the LW theorem in 
\cite{beuk_al}, one can 
only achieve that $u_n,v_n\in\Q$, and denominators have to be taken care of.

\section{Adapting Hilbert's proof to FPS}\label{sec_klazar}

In this section, we present another uncountable-sets-free 
proof of the transcendence of $\mathrm{e}$. It is due to this 
author. We adapt Hilbert's proof to FPS. We work in the domain $\R[[x]]$ 
of FPS with real coefficients and in the domain $\R\{x\}$ of those FPS in 
$\R[[x]]$ that are convergent. The latter are FPS with an infinite 
radius of convergence. For $f=f(x)\in\R[[x]]$ and $n\in\N_0$, 
we let $[x^n]f$ and $[x^n]f(x)$ denote the coefficient of $x^n$.

Let $f(x)=\sum_{n\ge0}a_nx^n\in\R[[x]]$ and $a\in\R$. We define the FPS
$$
{\textstyle
f(ax):=\sum_{n\ge0}a^na_nx^n\,.
}
$$
Then for every 
$a,b,c\in\R$ and $f(x),g(x)\in\R[[x]]$ we have 
$f(1x)=f(x)$, 
$$
(af+bg)(cx)=af(cx)+bg(cx),\,\text{ and }\,(fg)(cx)=f(cx)g(cx)\,. 
$$
Also, $f(0x)=[x^0]f$. If $f(x)\in \R\{x\}$ 
then $f(cx)\in\R\{x\}$. By 
$f(-x)$ we mean $f((-1)x)$.

In the 
next proposition, whose proof we omit, we 
collect properties of the 
unary operations of the derivative $f\mapsto f'$ and the primitive 
$f\mapsto\int f$. Recall that $\big(\sum_{n\ge0}a_nx^n\big)'=\sum_{n\ge0}(n+1)a_{n+1}x^n$ and 
$$
{\textstyle
\int\sum_{n\ge0}a_nx^n=\sum_{n\ge0}\frac{a_n}{n+1}x^{n+1}=
\sum_{n\ge1}\frac{a_{n-1}}{n}x^n\,.
}
$$

\begin{prop}\label{prop_formDeri}
Let $f,g\in\R[[x]]$ and $a,b\in\R$. The following holds.
\begin{enumerate}
\item $(af+bg)'=af'+bg'$, $(fg)'=f'g+fg'$ and $f(ax)'=af'(ax)$.
\item $\int(af+bg)=a\int f+b\int g$ and $\int f(ax)=\frac{1}{a}(\int f)(ax)$ if $a\ne0$.
\item  If $f\in\R\{x\}$ then $f',
\int f\in\R\{x\}$.
\item $\big(\int f\big)'=f$ and $\int f'=f-[x^0]f$.
\end{enumerate}
\end{prop}

We recall a~prominent FPS.

\begin{defi}\label{def_formExp}
The {\em FPS}
$${\textstyle
\mathrm{Exp}=
\mathrm{Exp}(x)=\sum_{n\ge0}\frac{1}{n!}x^n=
1+x+\frac{1}{2}x^2+\frac{1}{6}x^3+\ds\ \ (\in\Q[[x]])
}
$$
is called the (formal) exponential.
\end{defi}
It is easy to see that 
$\mathrm{Exp}\in\R\{x\}$, $\mathrm{Exp}'=\mathrm{Exp}$ and 
$\int\mathrm{Exp}=\mathrm{Exp}-1$.

We introduce a~family of 
{\em semiformal} unary operations on 
$\R\{x\}$. We call them semiformal because the coefficients of the 
FPS involved in these operations arise by limit transitions. In the 
usual formal operations with FPS, coefficients are obtained by finite expressions.

\begin{defi}[shifts]\label{def_formShif}
Let  $f(x)=\sum_{n\ge0}a_nx^n\in
\R\{x\}$ and let $b\in\R$. We define the {\em FPS} $f(x+b)\in\R\{x\}$ by
$$
{\textstyle
f(x+b):=\sum_{n\ge0}\big(\sum_{m\ge n}\binom{m}{n}a_mb^{m-n}\big)x^n\;.
}
$$    
\end{defi}
The coefficients in $f(x+b)$ are correctly defined because $f(x)$ is 
convergent. It is easy to show that $f(x+b)\in\R\{x\}$. Clearly, 
$f(x+0)=f(x)$ and if 
$f(x)\in\R[x]$ then $f(x+b)$ arises by the usual substitution in 
a~polynomial. We combine shifts and the operations $f(x)\mapsto f(ax)$ by 
means of brackets. For example, $f((x+a)+b)$ means $g(x+b)$ where 
$g(x)=f(x+a)$ or 
$f((ax)+b)$ means $g(x+b)$ where $g(x)=f(ax)$ or $f(a(x+b))$ means $g(ax)$ where $g(x)=f(x+b)$, and so on.
In the next two propositions, we show that shifts behave properly.  

\begin{prop}\label{prop_repeAdd}
Let $f(x)\in\R\{x\}$ and $c,d\in\R$. Then 
$$
f((x+c)+d)=f(x+(c+d))\,.
$$
    
\end{prop}
\proof
Let $f(x)=\sum_{n\ge0}a_nx^n$ and $f(x+c)=\sum_{n\ge0}b_nx^n$. We have 
$${\textstyle
b_n=\sum_{m\ge n}\binom{m}{n}a_m
c^{m-n}\,\text{ and }\,f((x+c)+d)=\sum_{n\ge0}\big(\sum_{m\ge n}\binom{m}{n}b_md^{m-n}\big)x^n\,.
}
$$
Then indeed
\begin{eqnarray*}
&&{\textstyle
[x^n]\,f((x+c)+d)=
\sum_{m\ge n}\binom{m}{n}\big(\sum_{l\ge m}\binom{l}{m}a_lc^{l-m}\big)d^{m-n}}\\  
&&={\textstyle
\sum_{m\ge n}\sum_{l\ge m}\frac{l!}{n!(m-n)!(l-m)!}a_lc^{l-m}d^{m-n}}\\
&&{\textstyle
=\sum_{l\ge n}\binom{l}{n}a_l\sum_{m=n}^l\binom{l-n}{l-m}c^{l-m}d^{m-n}}\\
&&{\textstyle
=\sum_{l\ge n}\binom{l}{n}a_l\sum_{j=0}^{l-n}\binom{l-n}{j}c^jd^{l-n-j}
}\\
&&{\textstyle=\sum_{l\ge n}\binom{l}{n}a_l(c+d)^{l-n}=[x^n]\,f(x+(c+d))\,.}
\end{eqnarray*}
The first equality follows from the definition of shifts. The second 
equality is an algebraic rearrangement. In the third 
equality, we change the order of summation. This is allowed because 
$f(x)$ is convergent. 
In the fourth 
equality, we introduce the variable $j=l-m$. The fifth equality follows 
from the binomial theorem. The last sixth equality follows from the 
definition of shifts.
\eproof

\begin{prop}\label{prop_AddaSoucin}
Let $f(x),g(x)\in\R\{x\}$ and $c\in\R$. Then 
$$
(f\cdot g)(x+c)=f(x+c)\cdot g(x+c)\;.
$$
    
\end{prop}
\proof
Let $f(x)=\sum_{n\ge0}a_nx^n$ and $g(x)=\sum_{n\ge0}b_nx^n$. We have  
$$
{\textstyle
[x^n]\,f(x+c)=\sum_{m\ge n}\binom{m}{n}a_mc^{m-n}
\,\text{ and }\,[x^n]\,g(x+c)=\sum_{l\ge n}
\binom{l}{n}b_lc^{l-n}\,.
}
$$
Then indeed 
\begin{eqnarray*}
&&{\textstyle
[x^n]\,f(x+c)g(x+c)=\sum_{k=0}^n\big(\sum_{m\ge k}\binom{m}{k}a_mc^{m-k}
\sum_{l\ge n-k}\binom{l}{n-k}b_lc^{l+k-n}\big)}\\
&&{\textstyle
=\sum_{k=0}^n\sum_{m\ge k}\sum_{l\ge n-k}\binom{m}{k}
\binom{l}{n-k}a_mb_lc^{l+m-n}}\\
&&{\textstyle
=\sum_{\substack{m,\,l\ge0\\m+l\ge n}}\sum_{k=0}^n
\binom{m}{k}\binom{l}{n-k}a_m b_lc^{l+m-n}}\\
&&{\textstyle
=\sum_{\substack{m,\,l\ge0\\m+l\ge n}}\binom{m+l}{n}a_mb_lc^{l+m-n}}\\
&&{\textstyle=
\sum_{j\ge n}\binom{j}{n}\big(\sum_{m=0}^j
a_mb_{j-m}\big)c^{j-n}=[x^n]\,(fg)(x+c)\,.}
\end{eqnarray*}
The first equality follows from the definitions of shifts and of multiplication of FPS. The second 
equality is an algebraic rearrangement. In the third 
equality, we change the order of summation. This is allowed because 
$f(x)$ and $g(x)$ are convergent. In the fourth 
equality, we use Vandermonde's identity (\cite{vdmonde}). In the 
fifth equality, we introduce the variable $j=l+m$. The last sixth 
equality follows from the 
definitions of shifts  and of multiplication of FPS.
\eproof

We prove a~curious semiformal exponential identity.

\begin{prop}\label{prop_expIdenForm}
Let $a,b\in\R$ with $a\ne0$. Then
$${\textstyle
\mathrm{Exp}\big((ax)+\frac{b}{a}\big)=\mathrm{e}^b\cdot
\mathrm{Exp}(ax)\,.
}
$$
\end{prop}
\proof
Indeed
\begin{eqnarray*}
{\textstyle
\mathrm{Exp}\big((ax)+\frac{b}{a}\big)}&=&{\textstyle \sum_{n\ge0}\big(\sum_{m\ge n}\binom{m}{n}(a^m/m!)(b/a)^{m-n}\big)x^n}\\
&=&{\textstyle
\sum_{n\ge0}\sum_{m\ge n}\frac{b^{m-n}}{(m-n)!}(a^n/n!)x^n}\\
&=&{\textstyle
\sum_{n\ge0}\big(\sum_{l\ge0}\frac{b^l}{l!}\big)
(a^n/n!)x^n}=\mathrm{e}^b\cdot
\mathrm{Exp}(ax)\,.
\end{eqnarray*} 
The first equality follows from the definitions of shifts and the 
operations $f(x)\mapsto f(ax)$.  The second equality is an algebraic 
rearrangement. In the 
third equality, we introduce the variable $l=m-n$. In the last
fourth equality, we take the constant $\mathrm{e}^b=
\sum_{l\ge0}\frac{b^l}{l!}$ out from the formal exponential. 
\eproof

\noindent
We use this identity only for $a=-1$.

We introduce a~family of semiformal functionals on 
$\R\{x\}$. 

\begin{defi}[Newton's $\int$]\label{def_NewtonInt}
Let $f(x)=\sum_{n\ge0}a_nx^n\in 
\R\{x\}$ and $u,v\in\R$. We define Newton's integral of 
$f(x)$ from $u$ to $v$ to be the difference
\begin{eqnarray*}
{\textstyle
\int_u^v f=\int_u^v f(x)}&:=&
{\textstyle\big(\int f\big)(v)-
\big(\int f\big)(u)
}\\&=&{\textstyle
\sum_{n\ge0}a_n\cdot\frac{v^{n+1}}{n+1}-\sum_{n\ge0}a_n\cdot\frac{u^{n+1}}{n+1}\,.
}
\end{eqnarray*}
\end{defi}
Since 
$\int f\in\R\{x\}$ (Proposition~\ref{prop_formDeri}), the integral is correctly defined. We obtain three transformations for it: linearity, additivity, and shift. The first two transformations are 
straightforward. 

\begin{prop}\label{prop_lineInte3}
Let $f,g\in\R\{x\}$ and $u,v,a,b\in\R$. Then
$${\textstyle
\int_u^v(af+bg)=a\int_u^v f
+b\int_u^v g\,.
}
$$
\end{prop}
\proof
This follows from the equality $\int(af+bg)=a\int f+b\int g$ in part~2 of Proposition~\ref{prop_formDeri}. 
\eproof

\begin{prop}\label{prop_addiInte3}
Let $f\in\R\{x\}$ and $u,v,w\in\R$. Then 
$${\textstyle
\int_u^w f=\int_u^v f+\int_v^w f\,.
}
$$
\end{prop}
\proof
Indeed, 
$$
{\textstyle
\big(\int f\big)(w)-\big(\int f\big)(u)=\big(\int f\big)(v)-\big(\int f\big)(u)+\big(\int f\big)(w)-\big(\int f\big)(v)\,. 
}
$$
\eproof

The shift transformation requires more work. 

\begin{prop}\label{prop_shifInte3}
Let $f\in\R\{x\}$ and $u,v,b\in\R$. Then
$${\textstyle
\int_u^v f(x+b)=\int_{u+b}^{v+b} f(x)\,.
}
$$
\end{prop}
\proof
Let $f(x)=\sum_{n\ge0}a_nx^n$. Then indeed
\begin{eqnarray*}
&&{\textstyle
\int_u^v f(x+b)}\\
&&{\textstyle
=\sum_{n\ge0}\big(\sum_{m\ge n}\binom{m}{n}a_mb^{m-n}\big)\frac{v^{n+1}}{n+1}-\sum_{n\ge0}\big(\sum_{m\ge n}\binom{m}{n}a_mb^{m-n}\big)\frac{u^{n+1}}{n+1}}\\
&&={\textstyle
\sum_{m\ge0}\frac{a_m}{m+1}\sum_{n=0}^m\binom{m+1}{n+1}b^{m-n}v^{n+1}-
\sum_{m\ge0}\frac{a_m}{m+1}\sum_{n=0}^m\binom{m+1}{n+1}b^{m-n}u^{n+1}}\\
&&{\textstyle
=\sum_{m\ge0}\frac{a_m}{m+1}\sum_{l=0}^{m+1}\binom{m+1}{l}b^{m+1-l}v^l-
\sum_{m\ge0}\frac{a_m}{m+1}\sum_{l=0}^{m+1}\binom{m+1}{l}b^{m+1-l}u^l}\\
&&{\textstyle
=\sum_{m\ge0}a_m\frac{(v+b)^{m+1}}{m+1}-\sum_{m\ge0}a_m\frac{(u+b)^{m+1}}{m+1}=\int_{u+b}^{v+b} f(x)\,.}
\end{eqnarray*}
The first equality follows from the definition of Newton's integral. In 
the second equality, we use the identity $\binom{m+1}{n+1}=
\frac{m+1}{n+1}\binom{m}{n}$ and change the order of summation. This is
allowed because the involved series absolutely converge. In the third 
equality, we introduce the variable $l=n+1$ and we add  
$$
{\textstyle
0=\sum_{m\ge0}\frac{a_m}{m+1}\binom{m+1}{0}b^{m+1}v^0-
\sum_{m\ge0}\frac{a_m}{m+1}\binom{m+1}{0}b^{m+1}u^0\,.
}
$$  
The fourth equality follows from the binomial theorem. The last fifth 
equality follows from the definition of Newton's integral. 
\eproof

We introduce another family of semiformal functionals on 
$\R\{x\}$.  

\begin{defi}[improper Newton's $\int$]\label{def_imprNewtonInt}
Let $u\in\R$ and $f(x)\in\R\{x\}$. If for every sequence $(b_n)\sus\R$ with $\lim b_n=+\infty$ the finite limit 
$${\textstyle
L=\lim_{n\to\infty}\int_u^{b_n}f\ \ (\in\R)
}
$$
exists, then $L$ does not depend on the sequence $(b_n)$ and we define 
$$
{\textstyle
\int_u^{+\infty}f=\int_u^{+\infty}f(x):=L\,.
}
$$
We call $\int_u^{+\infty}f$ improper Newton's integral of $f(x)$ from $u$ to $+\infty$.
\end{defi}
We leave the proof of
the independence of $L$ on $(b_n)$ for the 
interested reader as an exercise. 
We obtain the three transformations of 
linearity, additivity, and shift for the improper integral. All three are straightforward. 

\begin{prop}\label{prop_linInt4}
Let $f,g\in\R\{x\}$ and let $u,a,b\in\R$. Then the equality
$${\textstyle
\int_u^{+\infty}(af+bg)=a\int_u^{+\infty}f+b\int_u^{+\infty}g
}
$$
holds whenever the last two integrals exist.
\end{prop}
\proof
This follows from Proposition~\ref{prop_lineInte3} and 
the arithmetic of limits of real sequences. 
\eproof

\begin{prop}\label{prop_addiInte4}
Let $f\in\R\{x\}$ and $u,v\in\R$. 
The equality
$${\textstyle
\int_u^{+\infty}f=\int_u^v f+\int_v^{+\infty}f
}
$$
holds whenever one of the two improper integrals exists.
\end{prop}
\proof
Let $(b_n)\sus\R$ have $\lim b_n=+\infty$. Suppose that the former improper integral 
exists (the other case is similar). Thus the limit 
$L=\lim_{n\to\infty}\big(\int f\big)(b_n)$ ($\in\R$) exists. We have
\begin{eqnarray*}
&&{\textstyle
\int_u^{+\infty}f=
L-\big(\int f\big)(u)=
\big(\int f\big)(v)-\big(\int f\big)(u)+L-\big(\int f\big)(v)}\\
&&{\textstyle
=\int_u^v f+\int_v^{+\infty}f\,.
}
\end{eqnarray*}
The first equality follows from the definition of improper integrals. The 
second equality is an algebraic rearrangement. The last third equality 
follows from the definitions of proper and improper integrals.
\eproof

\begin{prop}\label{prop_shift4}
Let $f(x)\in\R\{x\}$ and $u,b\in\R$. The equality
$${\textstyle
\int_u^{+\infty}f(x+b)=
\int_{u+b}^{+\infty}f(x)
}
$$
holds whenever one of the two improper integrals exists.
\end{prop}
\proof
Let $(c_n)\sus\R$ have $\lim c_n=+\infty$.
Suppose that the former improper integral exists. Then also the latter improper integral exists and equals to it:
$$
\lim_{n\to\infty}{\textstyle
\int_{u+b}^{c_n}f(x)}\stackrel{\text{Prop.~\ref{prop_shifInte3}}}{=} 
\lim_{n\to\infty}
{\textstyle
\int_u^{c_n-b}f(x+b)=\int_u^{+\infty}f(x+b)\,.
}
$$
In the other way we go similarly:
$$
\lim_{n\to\infty}{\textstyle
\int_u^{c_n}f(x+b)}
\stackrel{\text{Prop.~\ref{prop_shifInte3}}}{=} 
\lim_{n\to\infty}
{\textstyle
\int_{u+b}^{c_n+b}f(x)=
\int_{u+b}^{+\infty}f(x)\,.
}
$$
\eproof

We need bounds on 
$\int_0^i p(x)\mathrm{Exp}(-x)$ 
for $i\in\N_0$ and $p(x)\in\Z[x]$. We obtain them from the next bounds.

\begin{prop}\label{prop_fakeLM}
Let $i,k\in\N_0$. Then
$${\textstyle
\big|\int_0^i x^k\cdot\mathrm{Exp}(-x)\big|\le i^{k+1}\mathrm{e}^i\,.
}
$$
\end{prop}
\proof
The integrand has the primitive  
$\sum_{n\ge0}\frac{(-1)^n}{n!(n+k+1)}x^{n+k+1}$. Thus 
$$
{\textstyle
\big|\int_0^i x^k\cdot\mathrm{Exp}(-x)\big|\le
\sum_{n\ge0}\frac{1}{n!}i^{n+k+1}+0=i^{k+1}\mathrm{e}^i\,.
}
$$
\eproof

We state and prove the version of Euler's identity for FPS. 

\begin{prop}\label{prop_formalEuId}
Let $k\in\N_0$. Then
$${\textstyle
\int_0^{+\infty} x^k\cdot\mathrm{Exp}(-x)=\int_0^{+\infty} \sum_{n\ge k}\frac{(-1)^{n-k}}{(n-k)!}x^n=k!\,.
}
$$
\end{prop}
\proof
We proceed by induction on $k$. Let 
$I_k$ denote the integral and let
$(a_n)\sus\R$ have $\lim a_n=+\infty$. From  
$$
{\textstyle
\int\mathrm{Exp}(-x)=\sum_{n\ge0}\frac{(-1)^n}{(n+1)!}x^{n+1}=1-\mathrm{Exp}(-x)
}
$$
we compute
$${\textstyle
I_0=\lim_{n\to\infty}
\big(1-\mathrm{Exp}(-a_n)-(1-\mathrm{Exp}(-0)\big)=1-\lim_{n\to\infty}
\mathrm{e}^{-a_n}=1\,.
}
$$
For $k\in\N_0$ we denote 
$F_k=x^k\cdot\mathrm{Exp}(-x)$.
Let $k\in\N$. We compute 
$$
F_k'=\big(x^k\cdot\mathrm{Exp}(-x)\big)'=kx^{k-1}\cdot\mathrm{Exp}(-x)-x^k\cdot\mathrm{Exp}(-x)=k\cdot F_{k-1}-F_k\,.
$$
The first equality follows from the definition of $F_k$. The second 
equality follows from part~1 of Proposition~\ref{prop_formDeri} and 
from the derivative 
$(\mathrm{Exp}(-x))'=-
\mathrm{Exp}(-x)$. The last third
equality follows from the definition of $F_k$. Hence 
$$
F_k=k\cdot F_{k-1}-F_k'\,. 
$$
Then indeed
\begin{eqnarray*}
I_k&=&\lim_{n\to\infty}{\textstyle
\int_0^{a_n}F_k}
=k\cdot\lim_{n\to\infty}{\textstyle\int_0^{a_n}F_{k-1}-}
\lim_{n\to\infty}{\textstyle\int_0^{a_n}F_k'}\\
&=&k\cdot I_{k-1}-\lim_{n\to\infty}\big(F_k(a_n)-F_k(0)\big)\\
&=&k\cdot(k-1)!-0=k!\;.
\end{eqnarray*}
The first equality follows from the definitions of $I_k$, $F_k$, and 
improper integrals. The second equality follows from the linearity 
of integrals and the identity $F_k=k\cdot F_{k-1}-F_k'$. In the 
third equality, we use the definitions of $I_k$, $F_k$, and 
improper integrals. We also use the definition of integrals and the fact 
that $\int F_k'=F_k-[x^0]F_k=F_k$ (part~4 of 
Proposition~\ref{prop_formDeri}). In the fourth equality, we use the 
induction, the value $F_k(0)=0$, 
and the limit $\lim_{n\to\infty}F_k(a_n)=0$. The last fifth equality is trivial.
\eproof

We conclude our FPS version of Hilbert's proof. We again assume for the contradiction that
$$
a_0+a_1\mathrm{e}+a_2\mathrm{e}^2+\ds+a_n\mathrm{e}^n=0
$$
for some $a_i\in\Z$, $n\in\N$, and $a_0\ne0$.  We multiply this equality 
by the integral 
$${\textstyle
\int_0^{+\infty}p_r(x)\cdot
\mathrm{Exp}(-x)\,,
}
$$
where $p_r(x)=x^r((x-1)(x-2)\ds(x-n))^{r+1}$ is as in 
Section~\ref{sec_intro}.  This integral exists by 
Propositions~\ref{prop_linInt4}
and \ref{prop_formalEuId}. Using  
Proposition~\ref{prop_addiInte4} we get for any $r\in\N$ the identity
\begin{eqnarray*}
0&=&A_r+B_r\\
&:=&{\textstyle\sum_{i=0}^n a_i\mathrm{e}^i\cdot\int_0^i p_r(x)\cdot\mathrm{Exp}(-x)+
\sum_{i=0}^n a_i\mathrm{e}^i\cdot\int_i^{+\infty} p_r(x)\cdot\mathrm{Exp}(-x)\,.
}
\end{eqnarray*}
The numbers $A_r,B_r\in\R$ are as in 
Section~\ref{sec_intro}, only now they are defined in an uncountable-sets-free 
way. The same contradiction as in Section~\ref{sec_intro} arises: every number
$B_r$ is an integral multiple of $r!$, $B_r\ne0$ for 
infinitely many $r\in\N$, and numbers $A_r$ are exponentially bounded. The 
computation of $B_r$ proceeds along the same lines as in 
Section~\ref{sec_intro}, only now we compute 
semiformally in $\R\{x\}$: 
\begin{eqnarray*}
B_r&=&{\textstyle
\sum_{i=0}^n a_i\mathrm{e}^i\int_i^{+\infty}p_r(x)
\cdot\mathrm{Exp}(-x)}\\
&\stackrel{\text{Prop.~\ref{prop_linInt4}}}{=}&{\textstyle
\sum_{i=0}^n a_i\int_i^{+\infty}p_r(x)\cdot
\mathrm{e}^i\cdot
\mathrm{Exp}(-x)}\\
&\stackrel{\text{Prop.~\ref{prop_expIdenForm}}}{=}&{\textstyle
\sum_{i=0}^n a_i\int_i^{+\infty}p_r(x)\cdot\mathrm{Exp}((-x)+(-i))}\\
&\stackrel{\text{Props.~\ref{prop_shift4} and \ref{prop_AddaSoucin}}}{=}&
{\textstyle
\sum_{i=0}^n a_i\int_0^{+\infty}p_r(x+i)\cdot\mathrm{Exp}(((-x)+(-i))+i)}\\
&\stackrel{\text{Prop.~\ref{prop_repeAdd}}}{=}&{\textstyle
\sum_{i=0}^n a_i\int_0^{+\infty}p_r(x+i)\cdot
\mathrm{Exp}(-x)}\,.
\end{eqnarray*}
Using Propositions~\ref{prop_linInt4} and \ref{prop_formalEuId}, and the definition of $p_r(x)$, we obtain the above properties of $B_r$.

It remains to bound, in an uncountable-sets-free 
way, the numbers
$${\textstyle
A_r=\sum_{i=0}^n a_i\mathrm{e}^i\int_0^i p_r(x)\cdot\mathrm{Exp}(-x)\,.
}
$$
The polynomial 
$$
p_r(x)=x^r((x-1)(x-2)\ds(x-n))^{r+1}
$$ 
has degree less than $(n+1)(r+1)$ and each coefficient in it is in absolute 
value at most $l^{r+1}$ where $l\in\N$ depends only on $n$.  By 
Propositions~\ref{prop_lineInte3} and \ref{prop_fakeLM}, each of the above  $n+1$ integrals in $A_r$ is in absolute value at most
$$
(n+1)(r+1)\cdot l^{r+1}\cdot n^{(n+1)(r+1)}\cdot\mathrm{e}^n\,.
$$
Thus $|A_r|\le c^r$ for every $r\in\N$ and some constant $c\ge1$. We get the same contradiction as 
in Section~\ref{sec_intro}.
\eproof

\end{document}